\documentclass{article}
\usepackage{graphicx}
\newtheorem{thm}{Theorem}

\newtheorem{lem}[thm]{Lemma}

\begin{document}
\title{Algorithms for recognizing knots and 3-manifolds}
\author{ Joel Hass \footnote{
This paper was completed while the author was visiting the Mathematical
Sciences Research Institute in Berkeley in 1996/7. Research at MSRI is
supported in part by NSF grant DMS-9022140.}}
\maketitle

\section{Algorithms and Classifications.}

Algorithms are of interest to geometric topologists for two reasons.
First, they have bearing on the decidability of a problem.  Certain topological
questions, such as finding a classification of four dimensional manifolds,
admit no solution.
It is important to know if other problems fall into this category.
Secondly, the discovery of a
reasonably efficient algorithm can lead to a computer program which can be
used to examine
interesting examples. In this paper we will survey some topological
algorithms, in particular
those that relate to distinguishing knots.  Our approach is somewhat
informal, with many details
omitted, but references are given to sources which develop these ideas in
full depth.

Given a question $Q$, a {\em decision procedure for $Q$} or 
an {\em algorithm to decide $Q$} can be
thought of as a computer program which will produce an answer to $Q$ in a
finite amount
of time.  A formal description of an algorithm or a computer is given by
the notion of a {\em Turing machine}.
A Turing machine is a basic computational device that reads and writes onto
a tape.
The questions such a
machine can decide are the same as those that can be decided by more
complicated computers.
The tape is divided
into cells, which the Turing machine can read from and write to, one at a
time. The tape has a
leftmost cell, but is infinite to the right. A finite set of symbols can be
written onto the tape - the usual English alphabet if we wish.  The Turing
machine has a finite number of possible states, and its behavior is
determined by its state. Initially, some finite number of cells on the tape
contain
symbols and the rest are blank.  At each time interval, the Turing machine
scans the
symbol at the current tape location, and in a manner determined by the
symbol and its
current state it: 1. Changes to a new state.  2.  Overwrites the symbol it
has read with
a new symbol.  3.  Moves the tape one cell left or right.  Some states are
final.  The computation
ends when they occur.

$Q$ is called {\em recursive} if there is an algorithm that produces an
answer in a finite amount of time. Showing that there is an algorithm to decide
$Q$ is equivalent to showing that $Q$ is recursive. It is often easier to
find an algorithm that
takes a finite amount of time to give a ``Yes" answer to a question, but
may run forever if the
answer is ``No".  This does not establish that the question is recursive. A
rather different
issue is the amount of time it takes an algorithm to answer a given
question.  This is
determined by its {\em computational complexity}.  Computational
complexity is a measure of the difficulty of a problem, and can have
implications to a real
world implementation of an algorithm, but is independent of the
decidability of a question.   A
good reference for the notions of Turing machine, algorithms and
decidability is \cite{HU}.

There are natural questions that do not admit any algorithm to decide them.
A famous
example is the word problem for finitely presented groups.  Given a
group described by a collection of generators and relations
$$G = \{ x_1, \dots , x_m ; r_1, \dots , r_n \},$$
this question asks whether a given product of
generators represents the identity element.  It was shown by Novikov and
Boone that there are groups in which there is no
algorithm to decide this question. An exposition can be found in \cite{Rotman}.
Closely related is the question of whether a finitely presented group $G$
is isomorphic to the
trivial group, which also cannot be decided by an algorithm.
Note however that it is easy to
construct an algorithm which will answer the triviality question with a
``Yes" in
finite time if the group is trivial. An algorithm which constructs all
products of
relations and their inverses and checks for the generators in this list will
accomplish this goal.  However this algorithm will run on forever if the
group is
non-trivial, so it does not decide the triviality question.

The notion of a classification is closely related.  Define a {\em
classification} of a set of objects to be a list containing each element of
the set once. 
Finding a
classification of some set of objects does not necessarily end its
mathematical interest.
As an example, it is easy to classify the natural numbers.

{\it I am indebted to W. Jaco and G. Kuperberg for helpful discussions.}

\section{Topological algorithms}

We will now consider the problem of classifying compact
orientable 3-manifolds. We seek a list containing each 3-manifold once, with the
property that if we are given a 3-manifold in some standard form then we can
determine where on the list it appears.

While surfaces have been classified for some time, a classification
remains elusive for 3-manifolds.  Markov showed that the Novikov-Boone
results implied that the
classification problem for 4-manifolds was not solvable \cite{Markov}.
Given a finitely
presented group, a compact 4-manifold (or $n$-manifold with $n>4$) can be
constructed with that
group as its fundamental group. Markov showed that a classification of
4-manifolds could be
used to give an algorithm to solve the problem of whether this presentation
defines
the trivial group, which we have seen is impossible.  While it seems
unlikely that a similar
type of problem arises for 3-manifolds, it is still unknown whether they
can be classified.

Here is an example of something that is not a classification of
3-manifolds.  All closed
3-manifolds can be triangulated \cite{Moise}.  Since there are only
finitely many ways
to glue together $k$ tetrahedra, we can construct all of these
systematically.  A
resulting complex is a manifold exactly when the links of all vertices are
2-spheres, a
condition that is easy to check.  Throwing away the non-manifolds gives a
method of
generating a list containing all 3-manifolds.  The drawback is that a given
manifold appears
many times in this list, and there is no known method to decide whether two
manifolds
in this list are homeomorphic.  However this procedure would lead to a
classification if
we could solve the \\
{\bf Recognition Problem for 3-manifolds:}  Give an
algorithm to decide whether two closed
3-manifolds are homeomorphic.  

The 3-manifolds are specified by the finite amount
of data needed to describe a finite triangulation.  This data can be
extracted from any
of the standard ways of describing 3-manifolds, such as surgery on a link,
Heegaard
diagrams, hierarchies, etc.  In dimension
three the PL category
is equivalent to the smooth or
topologically tame categories
\cite{Moise}.  For combinatorial manifolds, a classification is 
equivalent to a solution of the recognition problem.  Given a recognition algorithm, one 
can construct a list of manifolds using increasing numbers of simplices, discarding
duplicates by applying the recognition algorithm.  Conversely, a classification would 
give a recognition algorithm by comparing combinatorial manifolds to those in the list. While the
recognition problem in general is still open, important cases are known. 

\section{Surfaces in 3-manifolds}

A common approach in trying to understand 3-manifolds is to cut them open
along surfaces
into simpler building blocks, and to understand the ways that these are
recombined to
form the manifold. The cutting surfaces should reflect the global nature of the
3-manifold that they are dissecting, or the building blocks could become
more complicated than the
original manifold.    It appears to be counterproductive to
cut open along a surface with lots of knotted tubes, or with complicated
self-intersections, since the cut open 3-manifold would be more complex. To
get simplified pieces, several
possible cutting surfaces can be used. {\em Incompressible surfaces} are
embedded and
contain no trivial tubes or handles.  {\em Heegaard} surfaces, a second
important
class, cut the 3-manifold into two simple pieces, handlebodies.  Even with
these surfaces, careful choices
must be made for the cutting open procedure to be useful.

Once one has decided on a surface to cut along, there is still a great deal of
choice.  One could vary the surface in its isotopy class, perhaps creating
fingers which needlessly spiral around the 3-manifold.  It seems natural to
search for a
particularly simple representative in the surface's isotopy class. One
successful idea, developed
originally by Meeks and Yau, is to put a Riemannian metric on the
3-manifold and find a
surface of least area in the homotopy or isotopy class of the surface
\cite{MY}\cite{MSY}. It is a
non-trivial result that a least area surface tends to minimize its
self-intersections, as well as being rather rigidly situated in the
3-manifold.  In
Thurston's development of the theory of hyperbolic 3-manifolds, pleated
surfaces played
a similar role.  In the  piecewise linear (PL) context, where one has a
triangulated
3-manifold, an attempt to push the surface around until it becomes as
simple as possible
gives rise to what is called a {\em normal surface}. These ideas are
closely related -
normal surfaces are the discrete analogs of minimal surfaces.

A triangulated 3-manifold is a decomposition of a 3-manifold into a union
of tetrahedra,
which intersect one another along lower dimensional simplices.  We do not
restrict to a
combinatorial triangulation, so it's not forbidden for two tetrahedra to
intersect along
several faces or edges.  We can also generalize the triangulation to
allow ideal simplices, tetrahedra with some or all of their vertices removed. A
neighborhood of a vertex could then be a surface corresponding to a
boundary component
of $M$.

{\bf Definitions:}
{\em Normal triangles} are disks in a 3-simplex
which meet three edges and three faces of the 3-simplex, and {\em normal
quadrilaterals} are
disks in a 3-simplex which meet four edges and four faces of the 3-simplex.
An {\em elementary
disk} is a normal triangle or quadrilateral.  A {\em normal surface} in a
triangulated 3-manifold
is an embedded surface in $M$ which intersects each 3-simplex in a disjoint
union of elementary
disks.

\begin{figure}[hbtp]
\centering
\includegraphics[width=.6\textwidth]{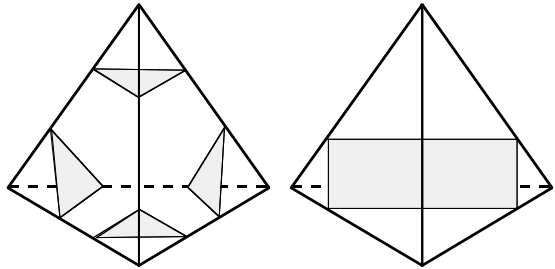}
\caption{Normal triangles and normal quadrilaterals.}
\label{normal}
\end{figure}

All four types of normal triangle can coexist disjointly in a 3-simplex.
However as soon
as one normal quadrilateral is around, the other two types of normal
quadrilateral can not
be present, or else an intersection would occur.

\begin{figure}[hbtp]
\centering
\includegraphics[width=.3\textwidth]{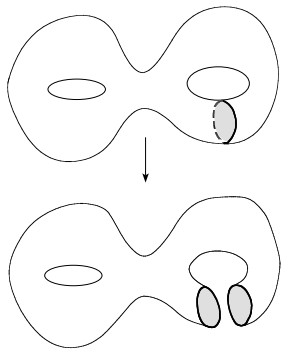}
\caption{A compressing disk for a surface and the result of a compression.}
\label{compress}
\end{figure}

\begin{figure}[hbtp]
\centering
\includegraphics[width=.8\textwidth]{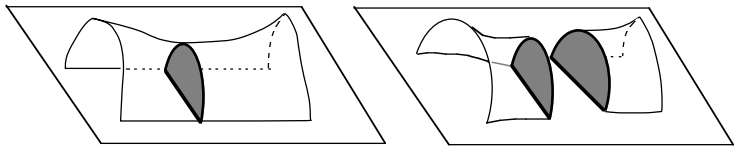}
\caption{A boundary compressing disk for a surface and the result of a
boundary compression.}
\label{boundary.compress}
\end{figure}

{\bf Definitions:}
A {\em compressing disk} for a surface $F$ inside a 3-manifold $M$ is an
embedded disk in $M$ which meets  $F$
along its boundary.  We
call the compressing disk non-trivial if the boundary curve of the disk
does not bound a disk on  $F$.  A surface
with no non-trivial compressing disks is {\em incompressible}.  A {\em
boundary compressing disk} for a surface
with boundary $F$ is an embedded disk in $M$ with an arc of its boundary on $F$
and the remainder of its boundary on
$\partial M$.  We call the boundary compressing disk non-trivial if the arc
on  $F$ is not parallel
to the boundary of $F$.  A surface with no non-trivial boundary compressing
disks is {\em
boundary incompressible}.  A non-trivial compressing disk can be used to
squeeze off a handle of
a surface, a process called {\em compression}, as in
Figure~\ref{compress}.  A similar process called {\em
boundary compression} squeezes arcs on a surface into the boundary of a
3-manifold, as in
Figure~\ref{boundary.compress}.

\section{Kneser's Theorem}

A 3-manifold $M$ which contains a separating 2-sphere can be cut open into two
3-manifolds with 2-sphere boundaries. By gluing in 3-balls to the boundaries, we get two
new closed 3-manifolds $M_1$ and $M_2$.  We say that the original manifold
is the
connect sum of these two manifolds, $M = M_1 \# M_2$, which are called the
summands.
This operation is well defined and associative if each manifold is oriented
and the gluing map of the
two 2-spheres is required to reverse orientation.  If the 2-sphere bounds a
ball, then one summand is
homeomorphic to
$M$ and the other to the 3-sphere.  This decomposition is called trivial.
The question arises
whether it's possible to keep on splitting $M$ indefinitely in a
non-trivial way, into more and
more pieces.

Normal surfaces were introduced by Kneser, who used them to prove the
following result
\cite{K}.

\begin{thm}
Let $M$ be a triangulated 3-manifold with $t$ 3-simplices and let $k(M) =
\mbox{dim}(H_1(M;Z_2)) + \mbox{dim}(H_1(M;Z)) + 6t$.
Them $M$ can be decomposed non-trivially along
2-spheres into at most $k(M)$ pieces.
\end{thm}

Suppose we have a collection of disjoint embedded 2-spheres in $M$.  We
will put the
2-spheres in a rigid position, making each of them a special type of surface.

In fact we will show that any embedded surface, not necessarily connected,
can be compressed and
isotoped to a union  (possibly empty) of normal surfaces. To do so we
introduce a notion of how
complex a surface is relative to a given triangulation.  The {\em weight}
$w(F)$ of a surface $F$ is the
number of times it intersects the 1-skeleton of $M$.  Weight gives an
analog of area
in the PL context.  It equals the area in the limiting case when all area
measure is concentrated near the 1-skeleton.

\begin{lem} \label{normalize}
Let $F$ be an embedded surface in $M$. Then after a series of compressions, isotopies and
removal of trivial 2-spheres,
$F$ becomes isotopic to a union (possibly empty) of disjoint normal surfaces.
\end{lem}
{\bf Proof:} Consider the intersection of $F$ with the triangulation. After a slight
perturbation of $F$, we can assume this intersection is transverse. We will simplify the intersection
by a process which we call {\em normalization}. 

$F$ intersects the interior of a tetrahedron $T$ in a collection of subsurfaces $F \cap T$, with boundary a
collection of disjoint simple closed curves. The boundary of $T$ is a 2-sphere, and each
simple closed curve in $F \cap \partial T$ cuts $\partial T$ into two disks. By applying the
Loop Theorem of Papakyriakopoulous \cite{Papakyriakopoulous} \cite{Hempel}, we can find a series of
compressions in $T$ which yield a new surface, all of whose components meet $T$ in disks and 2-spheres. 
Since $T$ is a ball, and therefore irreducible, the 2-sphere components are trivial and can be discarded. 
The weight does not increase in this process, though the number of components of $F$ may rise. Repeating for
each tetrahedra, we arrive at a surface $F_1$ of no higher weight with every curve in $F_1 \cap \partial T$
bounding a disk in $T$.

If a curve in $F_1 \cap \partial T$ lies completely inside a face of $T$, then we
can isotop the disk it bounds in $T$ across that face and eliminate the curve, as well as any
curves that lie inside it on the face. Repeating for other tetrahedra, we can assume that no such curves
exist in a face of any tetrahedron.

Now suppose that there is a curve $\delta$ of $F_1 \cap \partial T$ that meets an edge of $T$
in more than one point. If we consider all points of intersection of that edge with $F_1$,
we can find an adjacent pair of points which lies on one curve of $F_1 \cap \partial T$. 
These points can be connected by an arc on $F_1 \cap T$ whose interior is in the interior of $T$.
We can isotop this arc on
$F_1$ across the edge segment between these two points, reducing the weight of $F_1$ by two.

\begin{figure}[hbtp]
\centering
\includegraphics[width=.6\textwidth]{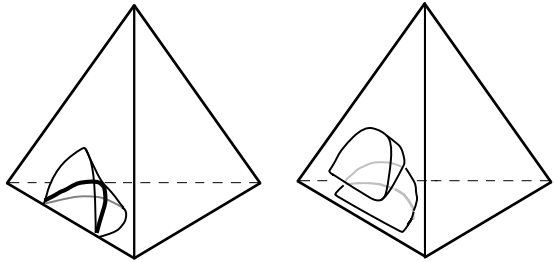}
\caption{An isotopy across an edge of a tetrahedron reduces weight.} 
\label{comp.edge} 
\end{figure}

Repeating in all the tetrahedra, we arrive at a surface  whose curves of intersection with the boundary of any
tetrahedron meet each edge at most once. A disk in a
tetrahedron whose boundary meets each edge at most once can only be a normal triangle or
quadrilateral. The resulting surface is a disjoint union of normal surfaces and surfaces
completely contained in a single tetrahedron. These latter surfaces can be compressed
to give trivial 2-spheres by an application of the
Loop Theorem, and the trivial 2-spheres can then be removed, leaving a normal surface as claimed. 

To illustrate the close connection between normal surfaces and minimal
surfaces, we now state an
important theorem of W. Meeks, L. Simon and S.T. Yau \cite{MSY}.
Lemma~\ref{normalize} is
essentially the same theorem in the PL context.

\begin{thm} \label{minimal}
Let $F$ be an embedded surface in a Riemannian 3-manifold $M$.  Then after
a series of compressions,
isotopies, and collapsing of the boundary of an I-bundle to its core, $F$
can be realized as a union
(possibly empty) of disjoint embedded minimal surfaces.
\end{thm}

The extra process of collapsing I-bundles can be seen in the Mobius band,
where a shortest curve
isotopic to the boundary can be homotoped to double cover the core.  In the
Riemannian setting, the metric
may force the minimizer to be a double cover.  In the PL setting, this type
of collapse is not
necessary, even where it is possible.

Each 3-simplex contains four distinct types of normal triangle and three
types of quadrilaterals.  A good way to keep track of
them is to notice that each normal triangle cuts off a unique vertex, and
each normal
quadrilateral separates one opposite pair of edges of the 3-simplex.   In a
given 3-simplex,
a union of normal triangles and quadrilaterals can contain at most five types,
- up to four triangles and at most one quadrilateral. Of course many
parallel copies of one
of the disk types can occur without causing intersections. The
complementary regions in the
3-simplex consist of a collection of product regions, \{triangle\}$\times$I or
\{quadrilateral\}$\times$I, together with at most six exceptional regions,
which are not
products.  The  exceptional regions meet either a vertex or at least two
distinct disk
types.  The two tetrahedra in Figure~\ref{normal} contain five and two
exceptional regions
respectively.

{\bf Proof of Kneser's Theorem:}  Let $ \{S_1, S_2, \dots, S_n\}$ be
disjoint 2-spheres in
$M$, no subset of which bounds a ball with some open balls removed, which
we call a {\em
punctured ball}. Let $F$ denote their union, and apply
Lemma~\ref{normalize} to $F$.  Then we
can isotope and compress to obtain a normal surface.  A compression causes
a 2-sphere $S$ to be
split into two 2-spheres $S_1$ and $S_2$.  The property that no subset of
the set of spheres
bounds a ball may not be preserved after the compression.  However if $S_1$
together with some
other spheres bounds a punctured ball $B_1$, and $S_2$ together with some
other spheres bounds a
punctured ball $B_2$, then $S$ together with some other spheres also bounds
a punctured ball.
So we can replace $S$ by one of $S_1$ and $S_2$ and get a new set of
spheres still having the
property that no subset bounds a punctured ball.  Repeating, we arrive  at
a collection of normal
2-spheres with this property which is as large as the initial collection.

We now show that if there are more than $k(M)$ disjoint normal 2-spheres,
then a subset
bounds a punctured ball, which implies there is a trivial $S^3$ summand in the
decomposition they define.  Each tetrahedron contains at most 6 non-product
regions.  The number of
components of $M - F$ which are not built entirely out of product regions
is bounded by
$6t$.  The product region components of $M - F$ are $I$-bundles with boundary a
2-sphere, and are either homeomorphic to
$S^2 \times I$ or to a non-trivial $I$-bundle over an embedded projective
plane in $M$.  Each of the
latter components contributes an $RP^3$ to the connect sum decomposition of
$M$, and a generator to
$H_1(M;Z_2)$. So the number of components of the complement of the
collection of 2-spheres is
bounded by $H_1(M;Z_2) + 6t$.  This gives a bound to the number of
separating 2-spheres in our
collection.  The number of non-parallel disjoint non-separating 2-spheres
is bounded by
$H_1(M;Z)$.  So if the number of 2-spheres is greater than
$k(M)$ then some subset of them must bound a punctured ball. 

Kneser's Theorem led to the establishment by Milnor of a unique
factorization theorem for 3-manifolds
into prime pieces.  See \cite{Hempel} for an exposition.

An almost identical argument proves a theorem of Haken  \cite{Haken2}.

\begin{thm}
Let $M$ be a triangulated 3-manifold with $t$ 3-simplices and let $k(M) =
\textstyle{dim}(H_1(M;Z_2)) + \textstyle{dim}(H_1(M;Z)) + 6t$. Then $M$
contains  at
most  $k(M)$ disjoint, non-parallel, incompressible surfaces.
\end{thm}

\section{Recognizing the unknot.}

Haken realized that the theory of normal surfaces had powerful applications
in the study
of 3-manifolds.  He used them in two important ways.
\begin{enumerate}
\item  To establish finiteness results about the number of ways in which
manifolds can be
cut open along surfaces other than 2-spheres. Haken showed that cutting
manifolds open
along suitably chosen {\em incompressible} surfaces resulted in a process
which terminated
after a finite number of steps \cite{Haken2}. Incompressible surfaces have
assumed a central role in the
theory of 3-manifolds.  The 3-manifolds which are irreducible and
which contain incompressible surfaces are called {\em Haken manifolds}.

\item To give algorithms to solve problems in 3-dimensional topology.

\end{enumerate}

We will concentrate on the second contribution, and describe
Haken's algorithm to recognize the unknot among the knots in the 3-sphere
\cite{H}.

A knot is the unknot if it bounds an embedded disk in $S^3$. The algorithm
will search for this
disk, and either produce it or show it does not exist in a finite amount of
time. In fact, the
same algorithm can find an embedded surface of smallest genus whose
boundary is the knot, giving the genus of
the know.  The knot is the unknot if and only if this smallest genus
surface is a disk. The surface will be
described as a normal surface.

To allow this, we need to extend the notion of normal surface to surfaces
and 3-manifolds with
boundary.  The definitions and pictures are the same, except that some
faces of some
tetrahedra lie on the boundary of the 3-manifold, giving a triangulation of
the boundary.
Lemma~\ref{normalize} generalizes with the extra operation of boundary
compression.

To get started, we need to describe a 3-manifold and a surface with a
finite amount of data.
A 3-manifold is described nicely by a triangulation.  The data is a finite
set of vertices, and
finite sets of pairs, triples and quadruples of vertices representing
edges, faces and
3-simplices.  There are some obvious conditions, such as that the three
pairs constructed from a
triple of vertices representing a face must occur as edges.  We allow these
sets of vertices to
contain a vertex or edge more than once, since our 3-simplices are not
required to be
combinatorial.

Describing a surface will need some additional data.  Taking advantage of
Lemma~\ref{normalize}, we work with normal surfaces.  For each tetrahedron
we have
seven types of normal triangles and quadrilaterals.  Assign to the $7t$
disk types found
in the 3-manifold the labels $\{\sigma_i\} $ and let $x_i$ denote
the multiplicity with which the disk type $\sigma_i$ occurs, $ \ 1 \le i
\le 7t$. The
vector of non-negative integers $x_i$ completely determines any embedded
normal surface.

We now ask which non-negative integer vectors $x_i$ give rise to a normal
surface.  Two conditions must be met.
\begin{enumerate}
\item Each tetrahedron contains at most one type of quadrilateral.  This is
called the
{\em quadrilateral} condition.
\item The disks must match up across a face separating two tetrahedra.
\end{enumerate}
The first condition means that if certain $\sigma_i$ have non-zero
coefficients, then others must have
zero coefficients.  The second means that the number of quadrilaterals and
triangles in a
tetrahedron having a given arc of intersection with a particular triangular
face must
equal the corresponding sum in the  tetrahedron adjacent across that face.
Since
each arc on a triangular face of a tetrahedron can be part of the boundary of
exactly one type of normal triangle and one type of normal quadrilateral,
this leads to
equations of the form
$ x_{i_1}  + x_{i_2} = x_{i_3} + x_{i_4}.$  These linear equations are
called the {\em
matching equations}. There are three of them for each pair of tetrahedral
faces which are glued to
one another.  Together with
the conditions $0 \le x_i$ we call them the {\em normal surface equations}.
There are no matching equations
associated to boundary faces.  Solutions of the normal surface equations
can have boundaries contained in these
faces.
If the quadrilateral condition and
the  normal surface equations hold, then the integers $x_i$ give a unique
embedded normal
surface, constructed by gluing together the normal triangles and
quadrilaterals in the unique way giving an embedding.  We denote by $X$ the
vector of
non-negative integers
$(x_1,x_2  \ldots, x_n)$, and also use $X$ to denote the embedded normal surface
made from the union of $x_i$ copies of the normal piece
$\sigma_i$.

\begin{figure}[hbtp]
\centering
\includegraphics[width=.6\textwidth]{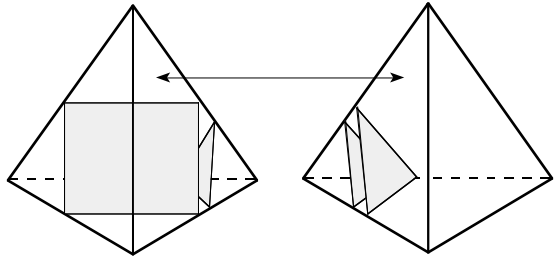}
\caption{Elementary disks in adjacent tetrahedra match up.}
\label{match}
\end{figure}

Two normal surfaces $A$ and $B$ may intersect one another. A surface $C$ is
obtained by taking their union and doing an operation called {\em regular
exchange}. This
involves cutting and pasting along double curves so as to preserve the
property that all
disks are normal. There is always a unique way to cut and paste two
normal disks intersecting along an arc to obtain a pair of disjoint normal
disks,
unless they are non-parallel quadrilaterals.  In that case no choice will
give normal disks.

\begin{figure}[hbtp]
\centering
\includegraphics[width=.6\textwidth]{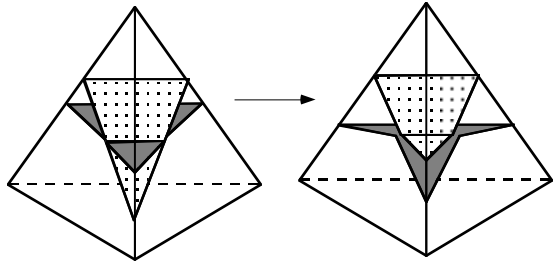}
\caption{A regular exchange keeps all disks normal.}
\label{regular}
\end{figure}

\begin{lem} \label{Euler}
If $A,B,C$ are solutions of the normal surface equations with $A+B=C$, and
$C$ gives rise to an embedded
normal surface, then so do $A$ and $B$.  Furthermore, $\chi (C) = \chi (A)
+ \chi (B)$ and
$w (C) = w(A) + w(B)$.
\end{lem}
{\bf Proof:} If $C$ corresponds to an embedded surface then it satisfies the
quadrilateral condition, so it doesn't contain intersecting quadrilaterals,
and neither do
$A$ or $B$.  We can calculate the Euler characteristic of $C$ by counting
vertices, edges and faces.  The
same number of each form $C$ and $A \cup B$.  The weight is unchanged when
regular
exchanges are made, and $C$ is obtained from $A \cup B$ by making regular
exchanges. 

A solution $C$ of the normal surface equations is {\em fundamental} if it
cannot be
written as a sum of two other solutions, $C \neq A + B$.

\begin{lem} \label{linear}
The set of fundamental solutions to the normal surface equations is finite.
\end{lem}
{\bf Proof:}  Real non-negative solution to the normal surface equations form a
cone contained in $R_+^n$. We can intersect this cone with the convex
simplex defined by
$\Sigma x_i = 1, x_i \ge 0$ and obtain a convex polyhedron with finitely
many vertices $
v_j$. The non-negative integer solutions are rational multiples of points
in this polyhedron.
Thus any solution of the normal surface equations can be
expressed as a rational linear combination of the  vertex solutions $ v_j
$.  Consider now
integral multiples $V_j = \lambda _j v_j $ of the $v_j$.  The normal
surfaces $V_j$ also form
a rational basis for all normal surfaces, so any normal surface $X$ can be
expressed in the
form $X = \Sigma t_j V_j $ with $0 \le t_j$.  Now consider the set
$S$ of real solutions to the normal surface equations which can be
expressed in the form  $s = \Sigma t_j V_j $ with $0 \le t_j \le 1, 1 \le j
\le n$.  $S$ is
compact, and so contains a finite number of integral points.   If
$X$ is an integral solution of the normal surface equations not in $S$,
then some $t_k > 1$ and we can
write
$X$ as the sum of two other positive integral solutions, $X = V_k + (\Sigma
t_j V_j - V_k)$.  Thus
all fundamental solutions lie in $S$, and there are only finitely many. 

Some subset of these fundamental solutions, those that also satisfy the
quadrilateral condition, correspond to
embedded surfaces.

The following result is due to Schubert
\cite{Schubert}.

\begin{lem} \label{schubert.sum}
If $C$ is a connected normal surface in an irreducible 3-manifold and
$C$ is not fundamental then we can find connected normal surfaces $A$ and
$B$ so that $C = A + B$.  If
$A$ and $B$ are chosen to minimize the number of curves of intersection in
$A \cap B$, then
no curve in $A \cap B$ is separating in both $A$ and $B$.
\end{lem}
{\bf Proof:}  Suppose $C = A_1 + A_2 \ldots + A_k$, where each $A_i$ is an
embedded normal
surface.  If $A_2$ intersects $A_1$, then we can perform regular exchanges
successively
along each curve in their intersection.  The number of surfaces adding to
$C$ changes by at
most one each time we perform a regular exchange.  Continuing for each
surface $A_j$, we eventually arrive at a
connected surface $C$, so at some point we must get two connected surfaces
whose sum is $C$. If these are
not embedded, we can perform regular exchanges along their curves of
self-intersection, giving a possibly
larger collection of embedded surfaces.  There are only finitely many
curves on which we can do a regular
exchange so the process stops with two embedded surfaces whose sum is
$C$.

Now pick $A$ and $B$ to minimize the number of curves of intersection in $A
\cap B$ among all pairs of connected
embedded normal surfaces whose sum gives $C$.  Suppose a curve $\alpha$ in
$A \cap B$ separates in
both $A$ and $B$.  Then a regular exchange along $\alpha$ cannot result in
a single connected surface.  It must
result in two connected surfaces $A'$ and $B'$, possibly non-embedded.
Regular exchanges on all
self-intersection curves of $A'$ and $B'$ results in embedded surfaces
$A''$ and $B''$, possibly not connected,
whose sum is $C$.  If $A'$ and $B'$ are not connected, then we can carry
out a series of regular exchanges
until we again arrive at exactly two connected surfaces. Eventually we
obtain $C$ as a sum of two embedded surfaces with fewer  curves of
intersection than $A$ and $B$, a
contradiction.  

A curve on the boundary torus of a knot complement is {\em essential} if it
does not bound a disk on the torus.
The following lemma allows us to reduce the search for an unknotting disk
to the finite collection of
fundamental surfaces.  Our proof is based on that of Jaco-Oertel\cite{JO}.

\begin{lem} \label{fundamental}
Suppose $C$ is a normal disk in a knot complement, with minimal weight among all
normal disks with essential boundary.  Then $C$ is fundamental.
\end{lem}
{\bf Proof:} If the lemma fails, then we can find two connected normal
surfaces $A$ and $B$ such that $C = A +
B$. Pick $A$ and $B$ as in Lemma~\ref{schubert.sum}. Since $\chi (C) = 1 =
\chi (A) + \chi (B)$, we have
several possibilities.   1.   $A$ is a punctured torus and $B$ a 2-sphere.
2. $A$ is a punctured Klein-bottle  and $B$ a 2-sphere.
3. $A$ is a disk and $B$ a torus.
4.  $A$ is a disk and $B$ an annulus.
5. $A$ is a disk and $B$ a Mobius band.
Possibilities involving
embedded  Klein bottles and projective
planes, though not more difficult, cannot occur in a knot complement.

In each case we can obtain a contradiction.  An argument used in \cite{JO} for
closed surfaces, extends to our setting.  This extension is explicitly derived in \cite{JT}. In each
of the above cases, if $C = A + B$ then we can find $A'$ and $B'$ such that
$C = A' + B'$ where
$A'$ is an essential disk and $B'$ is a torus or an annulus.  The essential normal disk $A'$
has lower weight than $C$, a contradiction.

Note: An alternate way to obtain a contradiction can be found by taking $C$ to be of minimal
complexity, as measured in the sense of Jaco-Rubinstein \cite{JR}.  In this setting the complexity
measures the length of intersection of a normal surface with the 2-skeleton of a triangulation.

\begin{thm} \label{unknot}
There is an algorithm to decide whether a knot is the unknot.
\end{thm}
{\bf Proof:} Triangulate the
knot complement $M_K$  and construct the finitely many fundamental
solutions. Among them find the ones
satisfying the quadrilateral conditions.  By calculating the Euler
characteristics, check if any
of these are disks.  If yes, test if the boundary of the disk is essential
on $\partial M_K$ by
checking whether it disconnects $\partial M_K$.  If there is a disk with
essential boundary on
$\partial M_K$ then $K$ is the unknot.  If there is no such disk it is not.

We now show why the algorithm works.  Suppose $K$ is the unknot.  Then it
bounds a disk. Lemma~\ref{normalize} implies that
there is a normal disk $D$ which bounds an essential curve on
$\partial M_K$.   Lemma~\ref{fundamental} implies we can find such a disk
among the
fundamental solutions.  So the algorithm finds an unknotting disk if it
exists. 

\section{Other algorithms to recognize the unknot.}

One might hope that an approach involving a search for moves that simplify
a knot
projection  would give an algorithm to recognize the unknot.  No approach
using this
idea has been found. The following  projection of the unknot, suggested by
Cameron
Gordon,  shows that if one wants to use a sequence of Reidemeister moves,
it may be
necessary to increase the number of crossings on the way to the unknot.  It
is an
interesting open problem to find a
larger class of allowable moves which would allow monotonic progress
towards the unknot.

\begin{figure}[hbtp]
\centering
\includegraphics[width=.5\textwidth]{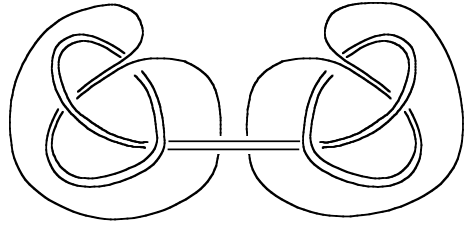}
\caption{The number of crossings in this projection of the unknot must be
increased before it can be decreased using Reidemeister moves.}
\label{crossings}
\end{figure}

Alternate algorithms to recognize the unknot could follow from a better
understanding of
certain knot invariants.  It is conjectured that the Jones polynomial of a
knot equals 1 if and only if the knot is the unknot.  Since there is a
simple finite
procedure for computing the Jones polynomial of a knot, a proof of this
conjecture would
provide an algorithm to decide if a knot is the unknot.

Kuperberg has pointed out that Thurston's work on geometric structures on
3-manifolds gives an alternate approach to constructing knot triviality
algorithms, one
of which we now describe.  We will refer to it as an {\em algebraic algorithm}.

This algorithm considers PL knots, as before. We triangulate the 3-sphere
and consider a
knot contained in the 1-skeleton.  The algorithm tests whether this knot is
the unknot by
running two processes in  parallel.

The first process checks if the
knot is unknotted by looking for embedded disks in the 2-skeleton.  If it
fails to find
any, it takes a barycentric subdivision and repeats.  The process ends if
it finds a
disk, in which case the knot is unknotted.

The second process searches for a non-cyclic finite representation of the
fundamental
group of the knot, in the following way.  It first computes the
Wirtinger presentation of the fundamental group of the knot complement.  It then
computes all homomorphisms of the group to the groups $S_3, S_4, S_5
\ldots$, where $S_n$
is the group of permutations of
$n$ elements. All homomorphisms to $S_n$  can be constructed by mapping the
Wirtinger
presentation generators to the finitely many possible sets of elements in
$S_n$, and
checking whether the relations of the knot group are satisfied in $S_n$.
The process
stops if it finds a homomorphism with non-cyclic image, in which case it
concludes
that the knot is non-trivial.

\begin{thm} \label{geometric.unknot}
The algebraic algorithm described above gives an algorithm to decide
whether a knot is
the unknot.
\end{thm}
{\bf Proof:}   If the knot is unknotted, the first process will end after a
finite amount of time.

If the knot is not the unknot, then it is a consequence of the existence of
geometric
structures on knot complements that the
fundamental group of the knot complement is a non-abelian residually finite
group
\cite{Hempel2}. {\em Residually finite} means that for any non-trivial element
of the
group, there is a homomorphism to some finite group which takes that
element to a
non-trivial element of the finite group.  There is no loss of generality in
considering as images only the groups $S_n$, which contain all other finite
groups
as subgroups.

A non-abelian residually finite group has a non-trivial commutator
subgroup.  An element
of this subgroup has non-trivial image in some $S_n$ under some
homomorphism. If the
image of all homomorphisms of the group were cyclic, and therefore abelian,
then this
element would always map to the trivial element, contradicting residual
finiteness.  So
for some element and some $S_n$ the image of the group is non-cyclic.

The process is guaranteed to stop after a finite time in either case.

\section{Classifying knots.}

To classify all knots we need to be able to decide not just whether a knot is
the unknot, but whether two arbitrary knots $K_1$ and $K_2$ are the same.
This was
carried out by Haken and Hemion using some additional arguments based on normal
surfaces \cite{He}.  We outline very briefly here an extension of the
algebraic algorithm
described above that can also give such a classification. This emerges from work
of Thurston, and was also related to me by Kuperberg.  We call it the {\em
geometric algorithm} to
distinguish knots.

Assume we are given two knots, $K_1$ and $K_2$, each presented as an
embedded circle in
the 1-skeleton of a triangulation of the 3-sphere. The idea is to set two
processes in motion, one of which will terminate if the knots are the same,
the other
if they are different.  Call the two triangulations $T_1$ and $T_2$.  The
first process
checks if $T_1$ is equal to $T_2$ by an isomorphism carrying $K_1$ to
$K_2$.  If not, it
performs a bistellar move on $T_1$ and checks again.  After trying all
bistellar moves on $T_1$, it
performs all pairs of bistellar moves.  Eventually it will get to all
sequences of $k$
successive bistellar moves.  If the two knots are the same, this will
terminate in
finite time.

The second process will terminate in finite time if the knots are distinct.
It has two
stages, first constructing geometric structures on each knot complement,
and then
checking if they are distinct.  It generates a series of subprocesses, run
in parallel,
as it tries different triangulations. It searches for the geometric
structures on each
knot complement, which must exist by Thurston's geometrization theorem.
The geometric
structure is constructed by triangulating the complement, with the
tetrahedra allowed to
be ideal, and constructing geometric structures on each tetrahedron so that
the angles
at edges and vertices match up.   If no compatible collection of angles can be
constructed, then the process looks for tori and annuli in the 2-skeleton
along which to cut up the knot
complement, and tries to geometrize each piece.  If it fails, it subdivides
barycentrically and repeats.  After a finite amount of time, geometric
structures on each knot
complement will be found by this process.  Moreover the meridian of the knot and the gluing
maps along splitting tori can be marked and remembered. It remains to check
if these structures
determine distinct knots.  This can be done by putting an
$\epsilon$ net on one 3-manifold, and trying to construct an isometry of
this net into the
second knot complement, preserving markings.  If this is impossible, then
the knots are
distinct.  If possible, then the process picks
$\epsilon$ smaller and tries again.  If the knots are distinct, the process
will fail to
construct an isometry eventually, and so will terminate in finite time,
establishing
that the knots are different.

\section{Other recognition results.}

In this section we briefly describe some other results on the problem of
recognizing knots
and three dimensional manifolds.

Schubert developed an alternate view of normal surfaces based on handle
theory.  He gave an algorithm to decide if a link is split \cite{Schubert}.
This approach is also used in
Jaco and Oertel
\cite{JO}, where an algorithm to decide if a 3-manifold is Haken is developed.
Haken and Hemion solved the recognition problem for Haken manifolds \cite{He}.
Rubinstein extended the notion of normal surfaces to the
concept of {\em almost normal surfaces}. He used this to describe an
algorithm which
decides whether a 3-manifold is homeomorphic to the 3-sphere
\cite{R2}.  Thompson gave an elegant argument that this algorithm works by
using the
notion of thin position
\cite{T}.
Rubinstein also described algorithms to recognize Lens spaces
and other 3-manifolds of small
genus.  See also Stocking
\cite{Stocking}. Rubinstein and Rannard have recently announced an
algorithm to recognize
Seifert Fibered manifolds.  Sela gave an algorithm to decide whether two
3-manifolds with
Gromov hyperbolic fundamental group are homotopy equivalent \cite{Sela}.
Birman and Hirsch \cite{BH}
have recently announced a new algorithm, based on work of Birman and
Menasco, which detects
whether a knot presented in braid form is the trivial knot.
Jaco and Tollefson describe algorithms to construct maximal families
of 2-spheres in a 3-manifold in \cite{JT}.  They also develop the important
idea of a {\em vertex surface}, a type of fundamental normal surface
introduced by Jaco and Oertel in \cite{JO},
which gives a more specialized representative
for a class of surfaces.
In \cite{HasLagPip} a bound for the complexity of the unknotting algorithm
is given.  It is also shown that this problem is in the class {\bf NP}.
Casson has recently announced a computation of the 
compexity of th 3-sphere recognition algorithm, which shows that it is
also in the class {\bf NP}.

\begin{flushleft}
Joel Hass\\ Department of Mathematics\\ University of California\\
Davis, CA 95616\\ e-mail: hass@@math.ucdavis.edu\\
\end{flushleft}


\begin{thebibliography}{HHH}

\bibitem{BH} J. Birman and M. Hirsch, Recognizing the Unknot, preprint.
\bibitem{H} W. Haken, Theorie der Normalflachen, Acta. Math. Vol. 105,
1961, 245-375.
1961, 245-375.
\bibitem{Haken2} W. Haken, Connections between topological and group
theoretical decision problems, Word
Problems, W.W. Boone (Ed.) North Holland, Amsterdam 1973, 427-441.
\bibitem{HasLagPip} J. Hass, J. Lagarias and N. Pippenger,
The computational complexity of knot and link problems, preliminary report, 
Proc. 38th Annual Symposium on Foundations of Computer Science, 1997, 172-1181.
\bibitem{He} G. Hemion, The Classification of Knots and 3-dimensional
Spaces, Oxford
University Press, 1992.
\bibitem{Hempel} J. Hempel, {\it 3-manifolds,} Annals of Math Studies 86,
Princeton U. Press 1976.
\bibitem{Hempel2} J. Hempel,
Residual finiteness for 3-manifolds,
{\it Combinatorial group theory and topology,} Alta, Utah, 1984,
Ann. of Math. Stud., 111, Princeton, 379-396.
Univ. Press, Princeton, NJ, 1987.
\bibitem{HU} J. Hopcroft and J. Ullman, Introduction to automata theory,
languages
and computation, Addison Wesley 1979.
\bibitem{JO} Jaco, W. and Oertel, U., An Algorithm to Decide if a
3-Manifold is a Haken
Manifold, Topology Vol. 23, No. 2, 1984, pp. 195-209.
\bibitem{JR} W. Jaco and J.H. Rubinstein, A piecewise linear theory of
minimal surfaces
in 3-manifolds,  J. Diff. Geom. Vol. 27, (1988), 493-524.
\bibitem{JT} W. Jaco and J.L. Tollefson, Algorithms for the
complete decomposition of a closed $3$-manifold, Illinois J. Math. 39
(1995), 358-406.
\bibitem{Johannson} K. Johannson, Classification problems in
low-dimensional topology,
Geometric and algebraic topology, Banach Center Publ., 18,
PWN, Warsaw (1986)  37-59.
\bibitem{K} H. Kneser,   Geschlossene Flachen in dreidimesionalen
Mannifgfaltigkeiten,
Jahresericht der Ent. Math. Verein 28 (1929) 248-260..
\bibitem{M}  S. V. Matveev, Algorithms for the recognition of the
three-dimensional sphere (after A. Thompson). Mat. Sb. 186 (1995) 69-84.
\bibitem{Markov} A.A. Markov, Insolubility of a problem of homeomorphy,
Proc. International Congress of
Mathematicians (1958) Cambridge Univ. Press, 300-306.
\bibitem{Moise} E.E. Moise, Affine Structures in 3-manifolds V, Annals of
Math. 55 (1952)
96-114.
\bibitem{MY} W. H. Meeks and S.T.Yau, Topology of three dimensional
manifolds and the embedding theorems in minimal surface theory, Annals of Math. 112
(1980) 441-484.
\bibitem{MSY} W. Meeks, L. Simon and S.T. Yau, Embedded minimal surfaces,
exotic spheres and manifolds of positive Ricci curvature, Annals of Math. 116 (1982) 621-659.
\bibitem{Papakyriakopoulous} C.D. Papakyriakopoulous, On Dehn's lemma and the asphericity of knots,
Annals of Math. 66 (1957), 1-26.
\bibitem{Rotman} J. J. Rotman, The Theory of Groups, Second Ed., Allyn and
Bacon (1973).
\bibitem{R1} J.H. Rubinstein, Polyhedral minimal surfaces, Heegaard
splittings and
decision problems for 3-dimensional manifolds,  Proceedings of the Georgia
Topology Conference, (1993) AMS/Intl. Press, 1-20.
\bibitem{R2} J.H. Rubinstein, An algorithm to recognize the 3-sphere, to
appear in
Proceedings of the ICM, Zurich 1994.  3.
\bibitem{Schubert} H. Schubert, Bestimmung der Primfaktorzerlegung von
Verkettungen, Math. Z. 76 (1961) 116-148.
\bibitem{Sela} Z. Sela, The isomorphism problem for hyperbolic groups I.
Annals of Math. 141 (1995) 217-283.
\bibitem{Stocking} M. Stocking, Almost normal surfaces in 3-manifolds.
Ph.D. Thesis, UC Davis (1996).
\bibitem{T} A. Thompson,  Thin Position and the Recognition Problem for
$S^3$ ,  Math. Research Letters 1, (1994) 613-630.
\bibitem{Thurston} W. Thurston, The geometry and topology of 3-dimensional
manifolds,
Princeton University Lecture Notes (1978).
\bibitem{Waldhausen} F. Waldhausen, Recent results on sufficiently large
3-manifolds,
 Proc. Sympos. in Pure Math. 32, Amer. Math. Soc. (1978) 21 - 38.
\bibitem{Waldhausen1} F. Waldhausen, Some problems on 3-manifolds,  Proc.
Sympos. in Pure
Math. 32, Amer. Math. Soc. (1978) 313-322.
\end{thebibliography}
\end{document}